\documentstyle[amscd,amssymb,verbatim]{amsart}

\newtheorem{theorem}{Theorem}[section]

\theoremstyle{definition}

\theoremstyle{remark}

\numberwithin{equation}{section}

\newcommand{\U}{U_q(G_2^{(1)})}
\renewcommand{\a}{\alpha}

\newcommand{\ep}{\epsilon}
\newcommand{\vep}{\varepsilon}

\begin{document}


\title[Quantum affine algebra $U_q(G_2^{(1)})$]
    {Level one representations of $U_q(G_2^{(1)})$
    }
\author{Naihuan Jing}
\address{Department of Mathematics\\
   North Carolina State University\\
   Raleigh, NC 27695-8205\\
   U.S.A.}
\email{jing\@eos.ncsu.edu}
\thanks{Research supported in part by NSA grants MDA904-96-1-0087
and MDA904-97-1-0062.}

\keywords{Quantum $Z$-operators, Vertex operators}
\subjclass{Primary: 17B}


\begin{abstract}
We construct a level one representation of the quantum affine
algebra $U_q(G_2^{(1)})$ by vertex operators from bosonic fields.
\end{abstract}

\maketitle

\section{Introduction} \label{S:intro}

Quantum affine algebras, the quantum groups associated to the affine
Kac-Moody Lie algebras, provide an important underlying
symmetry for the
quantum Yang-Baxter equation \cite{Dr} and quantum statistical
models \cite{JM}. Explicit realizations of their representations are much
needed in applications of quantum affine algebras. For instance, the
Frenkel-Reshetikhin vertex operators \cite{FR}
associated with the representations
can be used to give solutions of
the quantum Knizhnik-Zamoldchikov equation.

Lusztig first studied the abstract representations of quantum Kac-Moody
algebras \cite{L}. The program of constructing various representations was
started in
\cite{FJ} for level one irreducible modules of ADE types, and subsequently
twisted types were given in
\cite{J1} and $B_n^{(1)}$ in \cite{Br}. Recently we have constructed
symplectic quantum affine algebras in \cite{JKM2} for level one and
in \cite{JKM1} for level $-1/2$. The case of $F_4^{(1)}$
can also be done similarly \cite{J4} using the idea of quantum $Z$-algebras
\cite{J2, J4}. Besides the bosonic
constructions, fermionic
constructions were furnished in \cite{H}. The $q$-Wakimoto construction
was also known \cite{M, ABE, S} afterwards. Other representations
of classical quantum affine algebras have also been constructed
\cite{AOS}. The exceptional case of $G_2^{(1)}$ was the only case that has
not been explicitly constructed.

The purpose of the paper is to give a explicit level one
construction of the quantum
affine algebra $U_q(G_2^{(1)})$ by vertex operators.
The idea of the construction follows that of quantum $Z$-algebras \cite{J2,
J4},
which is a $q$-deformation of the classical ($q=1$) $Z$-algebras \cite{LW, LP}.
We construct some auxiliary vertex operators for the short root. This is
parallel to the known constructions of the affine Lie algebra $G_2^{(1)}$
\cite{BT, GNOS, XJ}, though the specialization of $q=1$ in our
construction is new even in the classical case.

The paper is organized as follows. In section two we review the quantum affine
algebra $U_q(G_2^{(1)})$. Section three gives the Fock space representation of
the quantum affine algebra $U_q(G_2^{(1)})$ stated in Theorem \ref{T:1}.
Section four uses quantum vertex operator techniques to prove Theorem
\ref{T:1}.
In the proof of Serre relations we have to show a relation about certain
symmetric functions, which is characteristic in the quantum affine algebras
as noted in \cite{J1}. The Serre relations in $G_2^{(1)}$
turn out to be the most
complicated one among both untwisted and twisted cases and actually
capture all the existed phenomena in other types.

\section{Quantum affine algebra $U_q(G_2^{(1)})$} \label{S:alg}

Let $\a_i$ ($i=1, 2$)
be the simple roots of the simple Lie algebra $G_{2}$, and $
\lambda_i$ be the fundamental
weight. Let $P={\bf Z}\vep_1+{\bf Z}\vep_2$ and
$Q={\bf Z}\a_1+{\bf Z}\a_2$ be the weight and root lattices.
We then let $\Lambda_i, i\in I=\{0, 1, 2\}$ be the fundamental
weights for the affine Lie algebra $G_2^{(1)}$, where $\Lambda_i=
\lambda_i+\Lambda_0$, and $\lambda_i$ are the fundamental weights
for the finite dimensional simple Lie algebra $G_2$.
The nondegenerate symmetric bilinear form $(\ |\ )$ on ${\bf h}^*$ is given by
\begin{equation}
(\alpha_i|\alpha_j)=d_ia_{ij}, \ \ (\delta|\alpha_i)=(\delta|\delta)=0
\ \ \mbox{for all} \ i,j
\end{equation}
where $(d_0, d_1, d_3)=(1, 1, 1/3)$ and $A=(a_{ij})=
\begin{pmatrix} 2 & -1\\
-3 & 2
\end{pmatrix}$.

Let $q_i=q^{d_i}=q^{\frac 12(\a_i|\a_i)}, i\in I$.
The quantum affine algebra $U_q(G_2^{(1)})$ is
the associative algebra with 1 over ${\bf C}(q^{1/6})$
generated by the elements $x_{ik}^{\pm}$, $a_{il}$, $K_i^{\pm 1}$,
$\gamma^{\pm 1/2}$, $q^{\pm d}$ $(i=1,2,\cdots,n, k\in {\bf Z},
l\in {\bf Z} \setminus \{0\})$ with the following defining relations
\cite{Dr, B, J3}:
\begin{align} \label{E:R1}
 [\gamma^{\pm 1/2}, u]&=0 \ \ \mbox{for all} \ u\in \textstyle {\bf U},\\
\mbox{} [a_{ik}, a_{jl}]&=\delta_{k+l,0}
\displaystyle\frac {[(\a_i|\a_j)k]}{k}
\displaystyle\frac {\gamma^k-\gamma^{-k}}{q-q^{-1}}, \label{E:R2}\\
\mbox{} [a_{ik}, K_j^{\pm 1}]&=[q^{\pm d}, K_j^{\pm 1}]=0, \label{E:R3}\\
 q^d x_{ik}^{\pm} q^{-d}&=q^k x_{ik}^{\pm }, \ \
q^d a_{il} q^{-d}=q^l a_{il}, \label{E:R4}\\
 K_i x_{jk}^{\pm} K_i^{-1}&=q^{\pm (\alpha_i|\alpha_j)} x_{jk} ^{\pm},
\label{E:R5}\\
\mbox{} [a_{ik}, x_{jl}^{\pm}]&=\pm \displaystyle\frac {[(\a_i|\a_j)k]}{k}
\gamma^{\mp |k|/2} x_{j,k+l}^{\pm},            \label{E:R6}
\\
(z-q^{\pm(\a_i|\a_j)}w)X_i^{\pm}(z)&X_j^{\pm}(w)
           +(w-q^{\pm(\a_i|\a_j)}z)X_j^{\pm}(w)X_i^{\pm}(z) =0, \label{E:R7}
\\
\mbox{}
\ [X_i^+(z),X_j^-(w)]
             =\frac{\delta_{ij}}{q_i-q_i^{-1}}&
                 \left(
                 \psi_i(w\gamma^{1/2})
                 \delta(\frac{w\gamma}{z})
                 -\varphi_i(w\gamma^{-1/2})
                  \delta(\frac{w\gamma^{-1}}z)
                 \right)                         \label{E:R8}
\end{align}
where
$X_i^{\pm}(z)=\sum_{n\in {\bf Z}}x_{i,n}z^{-n-1}$, $\psi_{im}$ and
$\varphi_{im}$ $(m\in {\bf Z}_{\ge 0})$
are defined by

\begin{align}
 \sum_{m=0}^{\infty} \psi_{im} z^{-m}
 &=K_i \textstyle {exp} \left( (q_i-q_i^{-1}) \sum_{k=1}^{\infty} a_{ik}
z^{-k}
\right), \\
 \sum_{m=0}^{\infty} \varphi_{i,-m} z^{m}
 &=K_i^{-1} \textstyle {exp} \left(- (q_i-q_i^{-1}) \sum_{k=1}^{\infty}
a_{i,-k} z^{k}
\right),
\\
 \sum_{r=0, \sigma\in S_m}^{m=1-A_{ij}}(-1)^r &
\left[\begin{array}{c} m\\r\end{array}\right]_i\sigma .X_i^{\pm}(z_1)\cdots
x^{\pm}_{i}(z_r)x^{\pm}_{j}(w)x^{\pm}_{i}(z_{r+1})\cdots x^{\pm}_{i}
(z_m)=0 . \label{E:R9}
\end{align}
where the symmetric group $S_m$ acts on  $z_i$
by permuting their indices.

\section{Fock space representations}   \label{S:fock}

Let $a_i(m)$ ($i=1, 2$) be the operators satisfying the Heisenberg relations
for $\U$ at
$\gamma=q$ and $b(m)$ and $c(m)$ be two independent free
bosonic operators with the relations:
\begin{align}
\ [a_i(m), a_j(n)]&=\delta_{m+l, 0}\frac{[(\a_i|\a_j)m]}{m}[m]\\
\ [b(m), b(n)]&=-\delta_{m+l, 0}\frac{[2m/3]}{m}[m]\\
\ [c(m), c(n)]&=\delta_{m+l, 0}\frac{[2m/3]}{m}[5m/3]\\
\ [a_i(m), b(n)]&=[a_i(m), c(n)]=[b(m), c(n)]=0 .
\end{align}
Let $\beta_2$ be an auxiliary simple root isomorphic to $\alpha_2$.
   We define the Fock module ${\mathcal F}$ as the tensor product of the
   symmetric
   algebra generated by $a_i(-n), b(-n), c(-n)$ ($n\in \mathbb N$) and
   the twisted group
   algebra $\mathbb C\{P+\mathbb Z\beta_2\}$ generated by $e^{\a}, e^{\beta}$
subject to
relation:
$$e^{\alpha_1}e^{\alpha_2}=-e^{\alpha_2}e^{\alpha_1}, \quad
e^{\alpha}e^{\beta}=e^{\beta}e^{\alpha}, \quad
e^{\alpha}e^{\alpha}=e^{\alpha}e^{\alpha}.
$$
where $\alpha\in P$, and $\beta$ is an element of the auxiliary lattice
$\mathbb Z\alpha_2$ (another copy of sublattice generated by
the short root $\alpha_2$). In the following we reserve $\beta$ to denote
an element from this auxiliary lattice.

The element
$1$ is the vacuum state. We define the action by
   $$ a_i(n).1 = 0 \quad (n>0)
      \; , \quad
      b_i(n).1 = 0 \quad (n>0) \;, $$
The element $a_i(0), b(0) $ act as differential operators  by
$$a_i(0)e^{\alpha}=(\alpha_i|\alpha)e^{\alpha}, \quad
b(0)e^{{\beta}}=-\frac23e^{{\beta}}.$$

  As usual we define the normal product as the ordered product by moving
annihilation operators $a_i(n), b(n), a_i(0), b(0)$ to the left.

   Let's introduce the following operators.
   \begin{align*}
   Y_1^{\pm}(z)&=
         \exp ( \pm \sum^{\infty}_{n=1}
                 \frac{a_1(-n)}{[n]} q^{\mp \frac{n}{2}} z^n)
         \exp ( \mp \sum^{\infty}_{n=1}
                 \frac{a_1(n)}{[n]} q^{\mp \frac{n}{2}} z^{-n})
         e^{\a_1} z^{\mp a_1(0)} , \\
   Y_2^{\pm}(z)&=
         \exp ( \pm \sum^{\infty}_{n=1}
                 \frac{a_1(-n)+b(-n)}{[n]} q^{\mp \frac{n}{2}} z^n)
         \exp ( \mp \sum^{\infty}_{n=1}
                 \frac{a_2(n)+b(n)}{[n]} q^{\mp \frac{n}{2}} z^{-n})\cdot\\
        &\qquad\qquad e^{\a +b} z^{\mp a_2(0\pm+b(0)} , \\
   U_{\pm}(z)&=
       \exp (\mp\sum^{\infty}_{n=1} \frac{[n/3]}{[2n/3]}b(\pm n) z^{\mp n})
       q^{\mp b(0)/2} , \\
   W_{\pm}(z)&=
       \exp (\mp\sum^{\infty}_{n=1} \frac{[n/3]}{[2n/3]}c(\pm n) z^{\mp n}).
   \end{align*}
   \bigskip
   \noindent
\begin{theorem} \label{T:1} The space
    ${\mathcal F}$ is a $\U$-module of level one
    under the action defined by
    $ \gamma \longmapsto q,
       K_i \longmapsto q^{a_i(0)}$,
    $ a_{im} \longmapsto a_i(m),
       q^d \longmapsto q^{\overline{d}}$, and
\begin{align*}
X_1^{\pm}(z) &\longmapsto Y_1^{\pm}(z)\\
X_2^{\pm}(z) &\longmapsto
      \frac{\pm Y^{\pm}_2(z)}{q_2-q_2^{-1}}
      \left(U_{\pm}(q^{\mp 5/6}z)W_{\pm}(q^{\mp 1/2}z)^{\pm 1}
        -U_{\pm}(q^{\pm 5/6}z)W_{\pm}(q^{\pm 1/2}z)^{\mp 1}\right).
\end{align*}
\end{theorem}

\section{Proof of the Theorem}

We now prove the theorem by checking that the action satisfy
Drinfeld relations. It is clear that the relations
(\ref{E:R1}-\ref{E:R5}) are true by the construction. The relation
(\ref{E:R6}) follows from the definition of $Y^{\pm}_i(z)$ and the
commutativity
among $\a_i(n)$, $b(n)$ and $c(n)$. So we only need to show the
relations (\ref{E:R7}-\ref{E:R8}).

We first compute the operator product expansions for $Y_i^{\pm}(z)$:
\begin{equation}\label{OPE1}
\begin{aligned}
Y_i^{\pm}(z)Y_j^{\pm}(w)&=:Y_i^{\pm}(z)Y_j^{\pm}(w):\\
   &\cdot  exp(-\sum_{n=1}^{\infty}\frac{[(\alpha_i|\alpha_j)n]}{n[n]}
     q^{\mp n}(\frac wz)^{n})z^{(\alpha_i|\alpha_j)}.\\
Y_i^{\pm}(z)Y_j^{\mp}(w)&=:Y_i^{\pm}(z)Y_j^{\mp}(w):\\
   &\cdot  exp(\sum_{n=1}^{\infty}\frac{[(\alpha_i|\alpha_j)n]}{n[n]}
     (\frac wz)^{n})z^{-(\alpha_i|\alpha_j)}.
\end{aligned}
\end{equation}

For $\epsilon=\pm=\pm 1$ we define
$$\begin{array}{rcl}
X^{+}_{2\epsilon}(z)&=&
   Y^{+}_2(z)U_{\epsilon}(q^{-5\ep/6}z)W_{\ep}(q^{-1\ep/2}z)\\
X^{-}_{2\epsilon}(z)&=&
 Y^{-}_2(z)U_{\epsilon}(q^{5\ep/6}z)W_{\ep}(q^{1\ep/2}z)^{-1}
\end{array}
$$
so that $X_2^{\pm}(z)=\frac1{q_2-q_2^{-1}} \left( X_{2+}^{\pm}(z)
-X_{2-}^{\pm}(z)
\right)$.

Note that for $i=j=1$ the relation (\ref{E:R7}) follows from the
$sl(2)$ case.
For $(\a_i|\a_j)=-1$ (i.e. $i\neq j$) equation (\ref{OPE1}) becomes
\begin{align}
Y_i^{\pm}(z)Y_j^{\pm}(w)&=:Y_i^{\pm}(z)Y_j^{\pm}(w):(z-q^{\mp 1}w)^{-1},\\
Y_i^{\pm}(z)Y_j^{\mp}(w)&=:Y_i^{\pm}(z)Y_j^{\mp}(w):(z-q^{\mp 1}w),
\end{align}
which implies that for $i\neq j$
\begin{align}
&(z-q^{\mp 1}w)X_i^{\pm}(z)X_j^{\pm}(w)=
(q^{\mp 1}z-w)X_j^{\pm}(w)X_i^{\pm}(z),\\
&[X_1^+(z), X_2^{-}(w)]=0,
\end{align}
where the latter is one case of relation (\ref{E:R8}).

To prove the remaining case of (\ref{E:R7}) we compute that
\begin{equation}\label{OPE2}
X_{2\ep}^+(z)X_{2\ep}^+(w)=:X_{2\ep}^+(z)X_{2\ep}^+(w):
\frac{z-w}{z-q^{2/3}w}q^{(1+\ep)/6}
\end{equation}
Then we immediately get the "+" case of relation (\ref{E:R7}) for $i=j=2$.
The "-" case is shown similarly.

In relation (\ref{E:R8}), again we only need to show the cases involved with
the short root $\a_2$, since the proof of $[X_1^+(z), X_1^{-}(w)]$ is
quite similar to that of type A \cite{J1}. Observe that
\begin{equation}\label{OPE3}
X^+_{2\ep}(z)X^-_{2,-\ep}(w)=:X^+_{2\ep}(z)X^-_{2,-\ep}(w):.
\end{equation}
Thus we reduce the relation to the commutators
$[X^+_{2\ep}(z), X^-_{2\ep}(w)]$.
We compute that
\begin{align*}
&[X^+_{2+}(z), X^-_{2+}(w)]\\
&=:X^+_{2+}(z)X^-_{2+}(w):
(\frac{z-q^{5/3}w}{z-qw}q^{-1/3}-\frac{w-q^{-5/3}z}{w-q^{-1}z}q^{1/3})\\
&=:X^+_{2+}(z)X^-_{2+}(w):\frac{z-q^{5/3}w}{z}q^{-1/3}\delta(\frac{qw}z)\\
&=(q^{-1/3}-q^{1/3})\psi_2(zq^{1/2})\delta(\frac{qw}z)
\end{align*}

Similarly we can prove that
\begin{equation*}
[X^+_{2-}(z), X^-_{2-}(w)]=(q^{1/3}-q^{-1/3})\phi_2(zq^{-1/2})
\delta(\frac{q^{-1}w}z).
\end{equation*}

Finally we use the quantum vertex operator calculus \cite{J1} to prove the
Serre relations.
The case ($a_{12}=-1$) is similar to that of $U_q(A_n^{(1)})$ \cite{J1}.
We only check the other one ($a_{ij}=-3$) in the ``+'' case:

\begin{multline}\label{Ser}
\sum_{\sigma\in S_4}\sigma .
\left(X_1^+(w)X^+_2(z_1)X^+_2(z_2)X^+_2(z_3)X_2(z_4)-\right.\\
[4]_2
X^+_2(z_1)X_1^+(w)X^+_2(z_2)X_2^+(z_3)X^+_2(z_4)+\\
 \frac{[4]_2[3]_2}{[2]_2}
X^+_2(z_1)X_2^+(z_2)X_1^+(w)X^+_2(z_3)X^+_2(z_4)-\\
[4]_2X^+_2(z_1)X^+_2(z_2)X_2^+(z_3)X_1^+(w)X^+_2(z_4)+\\
\left. X^+_2(z_1)X^+_2(z_2)X^+_2(z_3)X_2(z_4)X_1^+(w)\right)=0
\end{multline}

Recalling  $X_2^{\pm}(z)=\frac1{q_2-q_2^{-1}} \left( X_{2+}^{\pm}(z)
-X_{2-}^{\pm}(z)
\right)$
and using (\ref{OPE1}-\ref{OPE2}) and Wick's theorem we can reduce  the
left-hand side of (\ref{Ser}) into a linear combination of operator
product terms $:X_{1}^+(w)X^+_{2\ep_1}(z_1)X^+_{2\ep_2}(z_2)X^+_{2\ep_3}(z_3)
X_{2\ep_4}(z_4):$, where $\ep_i=\pm$. Thus the Serre relation is
equivalently reduced to four Serre-like relations grouped by the number of
appearance of
$X^+_{2+}(z_i)$ in the product. Due to (\ref{OPE3}) the most
complicated contraction functions comes from the case when all $\ep_i=+$.
All four subcases can be treated similarly. In the following we will only
prove the case when $\ep_i=+$.
Ignoring the factor $(q_2-q_2^{-1})^{-4}$ the left-hand side
 of this Serre-like
relation is  $q^2:X_{1}^+(w)X^+_{2+}(z_1)X^+_{2+}(z_2)X^+_{2+}(z_3)
X_{2+}(z_4):$ times the following expression.
\begin{align*}
&\sum_{\sigma\in  S_4}\sigma.\prod_i\frac1{(w-q^{-1}z_i)(z_i-q^{-1}w)}
\prod_{i< j}
\frac{z_i-z_j}{z_i-q^{2/3}z_j}\cdot\\
&\quad\left[(z_1-q^{-1}w)\cdots(z_4-q^{-1}w)
+[4]_2(w-q^{-1}z_1)(z_2-q^{-1}w)\cdots
(z_4-q^{-1}w)+\right.\\
&\quad\left[\begin{array}{c} 4\\2\end{array}\right]_2
(w-q^{-1}z_1)(w-q^{-1}z_2)(z_3-q^{-1}w)(z_4-q^{-1}w)+\\
&\left.\quad [4]_2(w-q^{-1}z_1)\cdots(w-q^{-1}z_3)
(z_4-q^{-1}w)+(w-q^{-1}z_1)\cdots (w-q^{-1}z_4)\right]
\end{align*}
where the symmetric group $ S_4$ acts on the ring of rational functions in
$z_i$ by permutations on the indices.
The $q$-binomial identity implies that the coefficients of $1$ and $w^4$ are
zero, and the expression in the bracket is then simplified to
\begin{align*}
&q_2^4(q_2^{-6}-1)\left[q_2^{-1}w^3\left(q_2^{-12}z_1-q_2^{-6}(1+q_2^{-2}+
q_2^{-4})z_2
+q_2^{-2}(1+q_2^{-2}+q_2^{-4})z_3-z_4\right)\right.\\
&\qquad +w^2(1+q^{-2}_2)\left(q_2^{-12}z_1z_2-q_2^{-6}(1+q^{-2}_2)z_1z_3+
q_2^{-4}(1+q^{-2}_2+q_2^{-4})z_1z_4+\right.\\
&\qquad \left. q^{-4}_2(1+q_2^{-2}+q_2^{-4})z_2z_3+q_2^{-4}(1+q_2^{-2})z_2z_4
-z_3z_4\right)+\\
&\qquad q_2^{-1}w
\left(q_2^{-12}z_1z_2z_3-q_2^{-6}(1+q_2^{-2}+q_2^{-4})z_1z_2z_3+\right.\\
&\qquad
\left.\left.q_2^{-2}(1+q_2^{-2}+q_2^{-4})z_1z_3z_4-z_2z_3z_4\right)\right]
\end{align*}
Let $f(z_1, z_2, z_3, z_4)$ denote the above expression. Since
$\prod_{i<j}(z_i-q_2^2z_j)(z_i-q_2^{-2}z_j)$ is symmetric, we see that
the Serre relation is equivalent to the following identity:
\begin{equation}\label{iden}
\sum_{\sigma\in  S_4}sgn(\sigma)\sigma \left(f(z_1, z_2, z_3, z_4)
\prod_{i<j}(z_i-q_2^{-2}z_j)\right)=0.
\end{equation}
We claim that (\ref{iden}) is true.  Notice that it is enough to check
 only the coefficients of $w$ and $w^2$, and even the two are quite similar.
The tedious checking
of the coefficient of $w$ shows that it is zero indeed.
Thus the Serre relation is proved.

The constructed
level one representation is reducible due to the presence of auxiliary
bosons $b(m)$ and $c(m)$. All integrable irreducible level one
modules are contained in the Fock representation and can be recovered by
the technique of the screening operators.


\end{document}